\documentstyle{article}
\newcounter{bean}

\def\X{{\cal X}}
\def\U{{\cal U}}
\def\A{{\cal A}}
\def\V{{\cal V}}
\def\F{{\cal F}}
\def\D{{\cal D}}
\def\O{{\cal O}}
\def\B{{\cal B}}
\def\M{{\cal M}}
\def\J{{\cal J}}
\def\H{{\cal H}}
\def\Z{{\bf Z}}
\def\H{{\bf H}}
\def\C{{\bf C}}
\def\HDI{{R}}
\def\R{{\bf R}}

\def\proclaim #1. #2\par{\medbreak{\bf#1.\enspace}{\it#2}\par
  \ifdim\lastskip<\medskipamount
  \removelastskip\penalty55\medskip\fi}
\newskip\Bigskipamount
\def\Bigskip{\vskip\Bigskipamount}
\begin{document}
\large
\centerline{\bf CUBICS, INTEGRABLE SYSTEMS, AND}
\bigskip
\centerline{\bf CALABI-YAU THREEFOLDS}

\vspace{0.3in}

\centerline{Ron Donagi\footnote{Supported by grants from NSF and
NSA}}
\vspace{0.1in}
\centerline{University of Pennsylvania}
\vspace{0.2in}
\centerline{and}
\vspace{0.2in}
\centerline{Eyal Markman\footnote{Supported by Rackham
Fellowship, University of Michigan, 1993}}
\vspace{0.1in}
\centerline{University of Michigan}

\vspace{0.3in}

\centerline{To Professor Hirzebruch, on his $65^{th}$ birthday}
\vspace{0.2in}
\S 0.  {\bf Introduction}

     In this work we construct an analytically completely
integrable Hamiltonian system which is canonically associated to
any family of Calabi-Yau threefolds.  The base of this system is
a moduli space of gauged Calabi-Yaus in the family, and the
fibers are Deligne cohomology groups (or intermediate Jacobians)
of the threefolds.  This system has several interesting
properties: \ \ the multivalued sections obtained as Abel-Jacobi
images, or ``normal functions'', of a family of curves on the
generic variety of the family, are always Lagrangian; the natural
affine coordinates on the base, which are used in the mirror
correspondence, arise as action variables for the integrable
system; and the Yukawa cubic, expressing the infinitesimal
variation of Hodge structure in the family, is essentially
equivalent to the symplectic structure on the total space.

     We begin our study by exploring this equivalence in a much
more general context.  The general question is as follows:\ \
Given a family of abelian varieties, when is there a symplectic
form on the total space with respect to which the abelian
varieties are Lagrangian?  This turns out to involve a symmetry
condition on the partials of the period map for the family.
Equivalently, this amounts to the existence of a field of cubics
on the tangent bundle to the base, such that the differential of
the period map at each point is given by contraction with the
cubic.  We give several versions of this `cubic condition' in
\S1.   For most of the examples and applications, we refer the
reader to the survey [DM].  The one example which we work out in
detail here is the integrable system of a Calabi-Yau family,
where the `cubic' is Yukawa's.

     The last section is entirely speculative.  We muse on the
role which the CY integrable system might play in the mirror
conjecture:  we would like to separate the conjecture into a
formal part, which deals with a `mirror transform' for objects
consisting of an integrable system together with certain
Lagrangian multisections and other lists of data, and a
geometric, Torelli-like part, saying that a CY family can be
recovered from its integrable system.  The attempt to read off
all relevant data (= the partition functions for the $A$ and $B$
models) directly from the system leads to some Hodge theoretic
questions which we cannot solve, concerning the information
encoded in a normal function and especially in its infinitesimal
invariant.

     It is our pleasure to acknowledge beneficial and pleasant
conversations with R. Bryant, I. Dolgachev, M. Green, P.
Griffiths, S. Katz, L. Katzarkov, D. Morrison,  T. Pantev and E. Witten.

\vspace{0.5in}

\S 1. {\bf The cubic condition}.

     Our goal in this section is to describe the conditions on a
given family $\pi: {\cal X} \to \U$ of abelian varieties (or
complex tori) for there to exist, locally over $\U$, an algebraic
(or analytic) symplectic structure $\sigma$ on ${\cal X}$ such
that $\pi$ is a Lagrangian fibration with respect to $\sigma$.
(In this situation we say that the family admits a
{\it Lagrangian structure}.)  What we want is a closed,
non-degenerate holomorphic $2$-form $\sigma$ on ${\cal X}$ such
that the fibers $X_u:= \pi^{-1}(u)$ are Lagrangian (= maximal
isotropic) subvarieties with respect to $\sigma$.

     By way of motivation, let us review the corresponding
questions for $C^{\infty}$ objects.  Any $C^\infty$ family $\pi:
\X \to \U$ of real tori is locally trivial, i.e. locally (over
$\U$) isomorphic to $\U \times (S^1)^n$.  In particular, a
Lagrangian structure exists whenever $\dim \U = n = \dim$
(fibers).
Given a $C^\infty$ symplectic manifold $(\X,\sigma)$ with
Lagrangian fibration $\pi: \X \to \U$, there is a natural affine
structure on the fibers $X_u$  [AG,2.42]; so if these fibers are
compact, they must be tori.  In this case the base $\U$ also
acquires a natural {\it integral-affine structure}, i.e. an
atlas of coordinate systems which differ on intersections by a
combination of arbitrary translations and ${\bf Z}$-linear
transformations.  To get these coordinates, we choose a basis
$\gamma_1(u),\cdots,\gamma_n(u)$ for the integral homology
$H_1(X_u,{\bf Z})$, varying smoothly with $u$ in an open subset
$\U^0$ of $\U$, as well as a $1$-form $\tau$ on $\pi^{-1}(\U^0)$
satisfying $d\tau = \sigma$.  (These always exist.) We then let
$$t_i(u):= \int_{\gamma_i(u)}\tau,\ \ \ \ \ i =
1,\cdots,n.\leqno(1)$$ The freedom of changing $\tau$ (by a
closed
$1$-form) and the $\gamma_i$ (by the action of $GL_n({\bf Z}))$
then allows the $t_i$ to change by a translation and
${\bf Z}$-linear transformation, respectively.  These $t_i$ are
called
{\it action variables}.  They define Hamiltonian vector
fields which are tangent to the $X_u$ and integrate to give the
{\it angle variables} $$(\varphi_1,\cdots,\varphi_n): X_u
\stackrel{\sim}{\rightarrow} {\bf R}^n/(2\pi{\bf Z})^n.$$ The
latter give the affine structure on the fibers.

     A {\it completely integrable Hamiltonian system}
(CIHS) on the $2n$-dimensional symplectic manifold $(\X,\sigma)$
is an $ODE$ (given, say, by a vector field $v$) which possesses
$n$ independent integrals $u_1,\cdots,u_n$ {\it in
involution} (i.e. they Poisson-commute).  This translates
immediately to a Lagrangian fibration $\pi: \X \to {\bf R}^n$
(where $\pi = (u_1,\cdots,u_n)$), and the flow is linearized on
the level tori $X_u$, i.e. $v$ is tangent to the fibers, and
$v|_{X_u}$ is a constant linear combination of the
$\partial/\partial\varphi_i$ with coefficients which depend
smoothly on $u$.

     The corresponding algebraic object is an
{\it algebraically completely integrable Hamiltonian
system} (ACIHS): we want $\X$ to be a $2n$-dimensional
non-singular complex algebraic variety, $\sigma$ a closed,
non-degenerate holomorphic $(2,0)$-form on $\X$, $\pi: \X \to \U$
a morphism whose fibers $X_u$ are abelian varieties which are
Lagrangian subvarieties of $(\X,\sigma)$, and $v$ a Hamiltonian
vector field for some algebraic (Hamiltonian) function on $\U$.
Such an $ODE$ can be solved more-or-less explicitly, in terms of
theta functions.  Numerous examples, both in mathematics and in
physics, are well-known; we refer to [DKN] for a survey, and to
[DM] for some recent examples.  What we are asking here is,
essentially, to decide when a given family of abelian varieties
underlies an ACIHS.

     The first new feature of the algebraic situation is that a
family of abelian varieties need not be locally trivial.  Rather,
it is specified by its {\it classifying map} $$q: \U \to {\cal
A}_g,$$ where $\A_g$ is the moduli space of abelian varieties of
dimension $g$ and of a given polarization type.  The given family
$\pi: \X \to \U$ is recovered as the pullback via $q$ of the
universal abelian variety $\X_g \to \A_g$.  (The latter, of
course, does not quite exist, but it exists locally away from
abelian varieties with excess automorphisms, and also globally on
a finite ``level'' branched cover of $\A_g$.)

     Let $\V \to \U$ be the {\it vertical} bundle determined by
$\pi: \X \to \U$.  Its fiber over $u \in \U$ can be identified
with the tangent space $T_0 X_u$, or with the universal cover of
$X_u$.  Deformation theory identifies the tangent space to $\A_g$
at a non-singular point $[X]\in \A_g$ with $Sym^2V$, where $V:=
T_0 X$.  We may thus write the differential of the classifying
map as $$dq: T\U \to Sym^2\V.\leqno(2)$$

     If the family $\pi: \X \to \U$ admits a Lagrangian structure
then the symplectic form $\sigma$ induces an isomorphism $$i =
\sigma^{-1}: \V^* \stackrel{\sim}{\rightarrow} T\U.\leqno(3)$$ We
refine our question by considering pairs consisting of the family
$\pi:  \X \to \U$ together with the isomorphism $i$, and asking
for the existence of $\sigma$ which is compatible with both.  Our
first answer is:
\proclaim{Theorem 1}. ({\bf Infinitesimal cubic condition})
Consider a
family $\pi: \X \to \U$ of polarized abelian varieties with
classifying map $q: \U \to \A_g$ and vertical bundle $\V$,
together with an isomorphism (3).  Then there exists a
non-degenerate holomorphic $2$-form $\sigma$ on $\X$ for which
$\pi$
is Lagrangian and which induces the given isomorphism $i$, if and
only if the composition of (2) and (3): $$dq \circ i \in
Hom(\V^*,Sym^2\V) = \Gamma(\V \otimes Sym^2\V)$$ comes from a
cubic $c \in \Gamma(Sym^3\V)$.  In this case, there is a unique
$\sigma$ for which the $0$-section of $\pi$ is Lagrangian.\par

\vspace{0.2in}

\noindent
{\bf Proof}.

     The short exact sequence of sheaves on $\X$: $$0 \to \pi^*\V
\to T\X \to \pi^*T\U \to 0$$ determines a subsheaf ${\cal F}$ of
$\Lambda^2T\X$ which fits in the exact sequences: $$0 \to \F \to
\Lambda^2T\X \to \pi^* \Lambda^2 T\U \to 0$$ $$0 \to
\pi^*\Lambda^2\V \to \F \to \pi^*(\V \otimes T\U)\to 0.$$ We are
looking for a $2$-vector $\sigma^{-1}$, a section of
$\Lambda^2T\X$, for which the fibers of $\pi$ are isotropic; this
means that $\sigma^{-1}$ goes to $0$ in $\pi^*\Lambda^2T\U$, so
it comes from a section of $\F$.  The question is therefore
whether the given $$i \in H^0(\U,\V \otimes T\U) \subset
H^0(\X,\pi^*(\V \otimes T\U))$$ is in the image of $H^0(\X,\F)$.
Locally in $\U$, this happens if and only if $i$ goes to $0$
under the coboundary map $$\begin{array}{ccc}\pi_*\pi^*(\V
\otimes T\U) & \rightarrow & \HDI^1\pi_*\pi^* \Lambda^2\V\\
\parallel & & \parallel\\ \V \otimes T\U & \rightarrow &
\Lambda^2\V \otimes \V.\\ \end{array}$$  The bottom map factors
as $\beta \circ (1 \otimes dq)$, where $$1 \otimes dq: \V \otimes
T\U \to \V \otimes Sym^2\V$$ comes from the differential of the
classifying map $q$, and $\beta$ is part of the (exact) Koszul
complex $$0 \to Sym^3\V \stackrel{\alpha}{\rightarrow} \V \otimes
Sym^2\V \stackrel{\beta}{\rightarrow} \Lambda^2\V \otimes \V \to
\cdots.$$ We conclude that the desired $2$-vector $\sigma^{-1}$
exists locally if and only if $$dq \circ i = (1 \otimes dq)(i)
\in \V \otimes Sym^2\V$$ comes, via $\alpha$, from a section of
$Sym^3\V$.  The global existence of $\sigma$ will follow by
patching local solutions, once we know uniqueness.  For this, let
$\sigma_1, \sigma_2$ be two local solutions.  Then
$\sigma_1 - \sigma_2$ induces the $0$-map from $T\U$ to $\V^*$,
hence it is the pullback of a $2$-form $\varphi$ on $\U$.  But
restricting to the $0$-section of $\pi$ shows $\varphi = 0$, so
$\sigma_1 = \sigma_2$ as required.
\begin{flushright}Q.E.D
\end{flushright}

     We view this theorem as an infinitesimal answer to our
problem:

     If we replace $\U$ by its first-order germ at the point $u
\in \U$, then the theorem provides a necessary and sufficient
condition for existence of a Lagrangian structure for $\pi: \X
\to \U$.  In general, it provides a non-degenerate $2$-form
$\sigma$, but an additional condition is needed for $\sigma$ to
be {\it closed}.  For this, it is convenient to replace the
moduli space $\A_g$ by the Siegel upper half space: $${\bf H}_g
:= \{g \times g\ \mbox{symmetric\ complex\ matrices\ \ Z \ with}
\ im(Z) > 0\}.$$ Over $\H_g$ there is (for each polarization
type)
a {\it universal} (marked) {\it abelian variety} $\X_g \to \H_g$,
whose fiber over $Z \in \H_g$ is $V/L$, where $V \approx {\bf
C}^g$ is a fixed complex vector space, and $L$ is the lattice
generated by the $2g$ columns of the $g \times 2g$ matrix
$(I,Z)$.  (Here $I$ is the $g \times g$ identity matrix, when the
polarization is principal, and in general it is the $g \times g$
integral diagonal matrix whose entries are the elementary
divisors of the polarization.)  Any family $\pi: \X \to \U$ of
abelian varieties which are {\it marked} (i.e. endowed with a
continuously varying symplectic basis of the fiber homologies
$H_1(X_u,\Z)$, $u \in \U$) determines a {\it period map} $$p: \U
\to \H_g,$$ and is recovered as the pullback via $p$ of the
universal family $\X_g$.  (The classifying map $q$, from $\U$ to
$\A_g$, is obtained by composition with the quotient map $\H_g
\to \A_g$ .)

     Let us consider the complex analogue of the action-angle
coordinates and the corresponding affine structures.  A
symplectic family $\pi: \X \to \U$ of bare abelian varieties,
given by the classifying map $q: \U \to \A_g$, determines locally
on $\U$ a set of $2g$ action variables, given by formula (1).  To
get action coordinates, we need to choose $g$ of these.  Now the
polarization on the fibers $X_u$ gives a symplectic structure on
$H_1(X_u,\Z)$, so the natural thing to do is to choose a
Lagrangian basis $\gamma_1(u),\cdots,\gamma_g(u)$, i.e. a basis
of a continuously varying Lagrangian subspace.  On the other
hand, the choice of a marking of the family $\pi$, i.e. a lifting
of $q$ to a period map $p: \U \to \H_g$, is equivalent to the
choice of a dual pair $\gamma_1,\cdots,\gamma_g$ and
$\gamma_{g+1},\cdots,\gamma_{2g}$ of Lagrangian bases.
Combining this with (1), we get:

\proclaim{Lemma 1}.  Let $\pi: \X \to \U$ be a family of abelian
varieties with Lagrangian structure.  The choice of Lagrangian
basis $\gamma_1,\cdots,\gamma_g$ of the fiber homologies
determines an affine structure on $\U$.  The choice of a marking
$\gamma_1,\cdots,\gamma_{2g}$ (or of a period map $p: \U \to
\H_g$) determines a pair of affine structures. The affine
coordinates $\{u_i\}$, $\{t_i\}$ corresponding to the Lagrangian
bases $\gamma_1,\cdots,\gamma_g$ and
$\gamma_{g+1},\cdots,\gamma_{2g}$, respectively, are related by:
$$dt_i = \sum p_{ij}(u)du_j.$$

     We note that, even with these choices, the angle coordinates
$\varphi_i$ are still multivalued, i.e. they can be defined only
locally.  Their differentials $d\varphi_i$ are more intrinsic:
they correspond symplectically to straight flows (with respect to
the affine structure) on $\U$.  The $\varphi_i$ can be uniquely
determined on the vertical bundle $\V$ (= the fiberwise universal
cover) by setting $\varphi_i = 0$ on the $0$-section of $\V$.

     We fix once and for all the vector space $V \approx {\bf
C}^g$.  Any marked family $\pi: \X \to \U$ of abelian varieties
then comes with a trivialization of its vertical bundle, $$\V
\approx V \otimes {\cal O}_\U.\leqno(4)$$ We identify $\H_g$ with
an open subset of $Sym^2V$, so
the differential of the period map becomes $$dp: T\U \to Sym^2V
\otimes \O_\U.$$ The affine structure on $\U$, determined by
$\gamma_1,\cdots,\gamma_g$ and a Lagrangian structure, is then
given by an isomorphism $$\alpha: V^* \otimes \O_\U
\stackrel{\sim}{\rightarrow} T\U.$$

\proclaim{Theorem 2}. ({\bf Global Cubic Condition})
\indent
Consider a family $\pi: \X \to \U$ of marked abelian
varieties, with period map $p: \U \to \H_g$, where $\U$ has an
affine structure $\alpha: V^* \otimes \O_\U
\stackrel{\sim}{\rightarrow} T\U$.  Then $\pi$ admits a
Lagrangian structure which induces $\alpha$ if and only if the
composition $$dp \circ \alpha \in Hom(V^* \otimes \O_\U,Sym^2 V
\otimes \O_\U) = \Gamma(V \otimes Sym^2V \otimes \O_\U)$$ comes
from a cubic $c \in \Gamma(Sym^3V \otimes \O_\U)$.\par
\eject

\noindent
{\bf Proof}.

     Note that the trivialization (4) takes $\alpha$ to the
isomorphism $i$ of (3), and $dp$ goes to $dq$ of (2).  The ``only
if'' therefore follows from Theorem 1.

     Conversely, we start with the standard symplectic structure
on the cotangent bundle $T^*\U$.  Via the affine structure
$\alpha$ and the trivialization (4), we get a symplectic structure
$\sigma$ on the vertical bundle $\V$.  The projection to $\U$ is
Lagrangian with respect to $\sigma$ as is the $0$-section.  We
want $\sigma$ to descend to $\X$.  Equivalently, $\sigma$ should
be invariant under translation by the section $$s_{m,n}: u
\mapsto m + p(u)n,$$ for any $m,n \in \Z^g$, which happens if and
only if the sections $s_{m,n}$ are Lagrangian.  Back on $T^*\U$,
we want the corresponding $1$-forms on $\U$ to be closed.  This
is always true for the $s_{m,0}$, since $\alpha$ is assumed to be
an affine structure.  For $s_{0,n}$, we need
closedness of the $g$ $1$-forms $\sum p_{ij}(u)du_j$.  This
amounts to equality of mixed partials, hence to the symmetry
condition on $dp \circ \alpha$.  We conclude that if $dp \circ
\alpha$ is symmetric then $\sigma$ descends to a Lagrangian
structure on $\pi$.  This induces $\alpha$ by construction.\\
\begin{flushright}Q.E.D
\end{flushright}

\vspace{0.2in}

\noindent
{\bf Remarks}.

     (1) Another way to state this is as follows.  Given the
period map $p$ and the isomorphism $i$, we get an isomorphism
$$\alpha: V^* \otimes \O_\U \stackrel{\sim}{\rightarrow} T\U.$$
If the cubic condition holds (with respect to $i$), we get the
$2$-form $\sigma$ by Theorem 1.  This $\sigma$ is closed (i.e.
gives a Lagrangian, or symplectic, structure) if and only if the
isomorphism $\alpha$ is integrable, i.e. if and only if $\alpha$
integrates to an affine structure on $\U$.

     (2) In the situation of the theorem, we deduce from the
equality of mixed partials that locally in $\U$ there is a
function $f \in \Gamma(\O_\U)$ such that $$\partial^2f = p,\ \ \
\ \partial^3f = c,$$ as sections of $Sym^2T^*\U$, $Sym^3T^*\U$
respectively.  Here $\partial$ denotes the total derivative with
respect to the affine structure on $\U$, so $\partial^k$ denotes,
essentially, the $k^{th}$ term in the Taylor expansion.  $f$ is
unique up to adding an affine function.  Conversely, any
holomorphic $f$ on a complex manifold $\U$ with $V^*$-affine
structure determines a Lagrangian family of abelian varieties on
the open subset of $\U$ where $Im(\partial^2f) > 0$.

     (3) Much of the above extends to Lagrangian structures on
families of complex tori.  The new feature is that $p$ lands in
$V \otimes V$, instead of $Sym^2V$, so we can write $p = p_+ +
p_-$, with $p_+$ symmetric, $p_-$ skew symmetric.  The result is
that in order to have a Lagrangian structure, $p_-$ must be
locally constant on $\U$, while $p_+$ satisfies a cubic condition
as above.

     In particular, everything goes through with no changes for
families of {\it polarized} complex tori, regardless of the
positivity of the polarization.  In terms of the period matrix,
we still require symmetry (Riemann's first bilinear relation) but
not positivity.  This is the situation which arises in \S2.

     (4) Another rather straightforward extension involves
replacing the symplectic structure $\sigma$ (or rather, its
inverse) by a Poisson structure, given by a $2$-vector $\psi \in
\Gamma(\Lambda^2T\X)$ satisfying a closedness condition.  Now we
start with an inclusion $i: \V^* \hookrightarrow T\U$ (or
$\alpha: V^* \otimes \O_\U \hookrightarrow T\U$) and we ask
whether a given family $\pi: \X \to \U$ of abelian varieties (or
complex tori) admits a Poisson structure $\psi$ such that
$\pi_*\psi = 0$ on $\U$ and which induces the given $i$.  As
before, the answer amounts to a cubic condition, $$dp \circ i \in
\Gamma (Sym^3\V).$$ We refer to [DM] for the details.

     (5) Dually, one can ask for existence of quasisymplectic
structures\\ (= closed, but possibly degenerate, $2$-forms) on
the
total space $\X$ of a family of abelian varieties, such that the
family is Lagrangian, but this time inducing a given inclusion
$$j: T\U \hookrightarrow \V^*.$$ Again, the answer [DM] is
in terms of a symmetry condition on $dp$, but of a slightly
different form:\ the condition is that $dp \otimes j$, which
lives in $$(Sym^2\V \otimes T^*\U) \otimes (\V^* \otimes
T^*\U),$$
maps under contraction to an element of $\V \otimes T^*\U \otimes
T^*\U$ which is symmetric in the last two variables, i.e. it is
in $\V \otimes Sym^2T^*\U$.  (This can be contracted with $j$
again, to give an element of $\otimes^3 T^*\U$.  The ``cubic
condition'', asserting that this last tensor is in $Sym^3T^*\U$,
follows from the previous symmetry condition, but is in general
too weak to imply it.)

     (6) Several examples and applications of the cubic
conditions are discussed in [DM] and [M2].  Jacobi's ACIHS (whose
flows include the geodesic flows on ellipsoids and the Neumann
flow) is particularly simple:\ the ``cubics'' turn out to be
Fermat cubics in $g$ variables, $\displaystyle{\sum^{g}_{i=1}}
x^3_i$.  These happen to have an invariance property which can be
stated as follows: the linear system of polar quadrics (generated
by the $g$ quadrics $x^3_i$) is invariant under the
$g$-dimensional group of rescalings, $x_i \mapsto a_i x_i$, $i =
1,\cdots,g$.  This in turn translates into existence of a
$g$-dimensional family of two-forms on the family for which the
fibration is Lagrangian, locally near each point of the base
$\U$.  (Only one of these, up to homothety, extends to a
Lagrangian structure on the entire family.)

     Jacobi's system has been extended to higher rank by Adams,
Harnad and Hurtubise [AHH] and by Beauville [B].  This system as
well as Hitchin's [H] are special cases of the spectral system
constructed in [M1] and [Bn] in terms of twisted endomorphisms of
vector bundles on curves.  The Poisson structure on this system
has a simple description in terms of the cubic condition, cf.
[DM].  A further generalization to a Lagrangian structure on the
space of Lagrangian sheaves on an arbitrary symplectic (or
Poisson, or quasi-symplectic) variety is given in [M2].  (This
includes the various spectral systems, as well as Mukai's
symplectic structure on spaces of sheaves on $K3$ surfaces,
etc.)  Here too, the symplectic (respectively Poisson,
quasi-symplectic) structure can be exhibited in terms of its
underlying
cubic, cf. [DM].

\vspace{0.3in}

\S 2. {\bf The Calabi-Yau integrable system}.

     In this section we construct the analytically completely
integrable Hamiltonian system associated to any complete family
of Calabi-Yau threefolds.  We also study the sections of this
system corresponding to curves on the generic Calabi-Yau.

     Let $X$ be a Calabi-Yau threefold, i.e. a three dimensional
compact K\"{a}hler manifold satisfying $$H^0(\Omega^1_X) =
H^0(\Omega^2_X) = 0$$ $$ \omega_X \approx \O_X.$$ Here
$\Omega^i_X$ is the sheaf of holomorphic $i$-forms on $X$, $\O_X
= \Omega^0_X$ is the structure sheaf, and $\omega_X = \Omega^3_X$
the canonical sheaf.  We denote the holomorphic tangent bundle by
$T_X$.  By (holomorphic) gauge we mean a non-zero section $s \in
H^0(\omega_X)$.

     Hodge theory gives us the Hodge decomposition $$H^3(X,{\bf
C}) = H^{3,0} \oplus H^{2,1} \oplus H^{1,2} \oplus H^{0,3}$$ and
the Hodge filtration $$F^p = \oplus_{i \geq p} H^{i,3-i},$$
together with bilinear maps $${\cal B}_p : H^1(T_X) \times
H^{p,q} \to H^{p-1,q+1}$$ which express the variation of the
Hodge structure: as $X$ is deformed in a direction $v \in
H^1(T_X)$, its $H^{p,q}$ moves to first order within $F^{p-1}$
(Griffiths transversality) and an element $w \in H^{p,q}$ moves
in the direction of ${\cal B}_p(v,w)$.  The $H^{p,q}$ can be
described as spaces of harmonic $(p,q)$-forms on $X$ (thus giving
the inclusion $H^{p,q} \hookrightarrow H^3(X,{\bf C}))$, or as
the sheaf cohomology groups $H^q(X,\Omega^p_X)$.  In the case of
a CY threefold $X$, we have:
\begin{list}{{\rm(\roman{bean})}}{\usecounter{bean}}
\item $H^{3,0} = H^0(\omega_{X})$ is the one-dimensional complex
vector space of holomorphic gauges.
\item $H^{2,1} = H^1(\Omega^2)$ is isomorphic to $H^1(T_X)$.
More canonically, each holomorphic gauge $s \in H^{3,0}$, $s \neq
0$, determines an isomorphism \begin{eqnarray*}\rfloor s:
H^1(T_X) &
\stackrel{\sim}{\rightarrow} & H^{2,1}\\v & \mapsto &
\B_3(v,s).\end{eqnarray*}
\item $H^{1,2}$ and $H^{0,3}$ can be identified either with the
complex conjugates of $H^{2,1}$ and $H^{3,0}$ (with respect to
the real subspace $H^3(X,\R)$ of $H^3(X,\C))$, or with their dual
spaces $(H^{2,1})^*$, $(H^{3,0})^*$ (the pairing is given by
integration over $X$).
\end{list}

     The Griffiths intermediate Jacobian of $X$ is
$$J(X):= H^3(X,\C)/(F^2H^3 + H^3(X,\Z)) \approx
(F^2H^3)^*/\lambda(H_3(X,\Z)),$$ where $H_3(X,\Z)$ is identified
with a lattice in $(F^2H^3)^*$ by the map $\lambda: H_3(X,\Z)
\hookrightarrow (F^2H^3)^*$ which is the composition of
integration $$\begin{array}{ccc}H_3(X,\Z) & \rightarrow &
H^3(X,\C)^*\\ \gamma & \mapsto & \int_\gamma\\ \end{array}$$
and projection $$H^3(X,\C)^* \rightarrow\!\!\!\rightarrow
(F^2H^3)^*.$$
$J(X)$ is a complex torus, but in general {\it not} an
abelian variety: it satisfies Riemann's first bilinear relation,
but not the second.  It has the important property that it
depends holomorphically on $X$.  This means that given a
holomorphic family $\chi: \X \to \M$ of compact K\"{a}hler
manifolds $X_t$, $t \in \M$, the intermediate Jacobians
$J(X_t)$ fit together naturally to form a holomorphic family $\pi:
\J \to \M$, called the relative intermediate Jacobian of the
family $\chi$.

     The theorem of Bogomolov, Tian and Todorov [Bo,Ti,To] says
that each moduli space $M$ of Calabi-Yaus is smooth and
unobstructed, i.e. the Kodaira-Spencer map gives an isomorphism
$$T_{[X]} M \stackrel{\approx}{\rightarrow} H^1(T_X)$$ for every
CY $X$.  We say that a family $\chi: \X \to \M$ of CYs $X_t,
t \in \M$, is {\it complete} if its classifying map $\M \to M$ is
a local isomorphism.  It follows that $\M$ is smooth with tangent
space $H^1(T_{X_t})$ at $t$.  Typically such families might
consist of all CYs in some open subset of $M$, together with
some ``level'' structure.  We fatten the family by adding the
holomorphic gauge: let $\rho: \widetilde{\M} \to \M$ be the
$\C^*$-bundle over $\M$ underlying the line bundle $\chi_*
\omega_{\X/\M}$, so a point of $\widetilde{\M}$ is given by a
pair $(X,s)$ with $X \in \M$ and $s \in H^0(\omega_X) \backslash
(0)$.  There is a pullback family $\widetilde{\chi}:
\widetilde{\X} \to \widetilde{\M}$, where $\widetilde{\X}:= \X
\times_\M \widetilde{\M}$.  We let $\pi: \J \to \M$ and
$\tilde{\pi}: \widetilde{\J} \to \widetilde{\M}$ denote the
relative intermediate Jacobians of $\chi,\widetilde{\chi}$
respectively.

\proclaim{Theorem 3}. $\pi: \widetilde{\J} \to \widetilde{\M}$ is
an analytically completely integrable Hamiltonian system.\par

\vspace{0.2in}

\noindent
{\bf Proof}.

     \underline{Step I}.  We claim that there is a canonical
isomorphism $$i^{-1}: T_{(X,s)} \widetilde{\M} \approx
F^2H^3(X,\C).\leqno(5)$$ Indeed, $\rho: \widetilde{\M} \to \M$
gives a short exact sequence: $$0 \to
T_{(X,s)}(\widetilde{\M}/\M) \to T_{(X,s)} \widetilde{\M} \to T_X
\M \to 0,$$ in which the subspace can be identified with
$H^0(\omega_X)$, and the quotient with $H^1(T_X)$, hence via
$\rfloor s$
with $H^1(\Omega^2_X)$.  We need to identify this sequence
naturally with the one defining $F^2 H^3$: $$0 \to H^0(\omega_X)
\to F^2H^3(X,\C) \to H^1(\Omega^2_X) \to 0.$$ In other words, we
are claiming that the extension data in the two sequences match,
globally over $\widetilde{\M}$.

     For this, it suffices to construct a natural map from
$T_{(X,s)} \widetilde{\M}$ to $F^2H^3(X,\C)$, which induces the
identity on the subspace and quotient.

     One way to obtain this map is to note that the bundle $\HDI^3
\widetilde{\chi}_*\C$ has the {\it tautological section} $s$,
which is actually in the subbundle $F^3\HDI^3
\widetilde{\chi}_*\C$.  The Gauss-Manin connection applied to
this section maps $T_{(X,s)} \widetilde{\M}$ to $H^3(X,\C)$, but
by Griffiths transversability it lands in $F^2H^3$.  This map is
easily seen to have the required properties.  For future use we
want an explicit description of this map.  We do this using
Dolbeault cohomology.

     Think of a $1$-parameter family $(X_t,s_t) \in
\widetilde{\M}$, depending on the parameter $t$, as a fixed
$C^\infty$ manifold $X$ on which there is given a $1$-parameter
family of complex structures, specified by their
$\overline{\partial}$-operators $\overline{\partial}_t$, and also
a family $s_t$ of $C^\infty$ $3$-forms on $X$ such that $s_t$ is
of type $(3,0)$ with respect to $\overline{\partial}_t$.  Since
the $s_t$ are now on a fixed underlying $X$, we can expand with
respect to $t$: $$s_t = s_0 + ta\ \ \ \ (\mbox{mod.}\ \
t^2).$$
Griffiths' transversality now says that the derivative, $a$, is
in $F^2H^3$.  It clearly depends only on first-order data, i.e.
on
the tangent vector to $\widetilde{\M}$ along $(X_t,s_t)$ at $t =
0$.  This produces the map $$T_{(X,s)} \widetilde{\M} \to
F^2H^3(X,\C)$$ with the desired properties.

\vspace{0.2in}

     \underline{Step II}.  We want a non-degenerate $2$-form
$\sigma$ on $\widetilde{\J}$ making $\pi$ Lagrangian.  By Theorem
1, this is equivalent to finding a cubic $\tilde{c}$ on
$T\widetilde{\M}$ satisfying $\rfloor\tilde{c} = d\tilde{p} \circ
i$,
where
$\tilde{p}$ is the period map for $\widetilde{\chi}:
\widetilde{\X} \to \widetilde{\M}$ and $i$ is the inverse of
(5).  Since $\tilde{p}$ factors through the period map $p$ for
$\chi: \X \to \M$, it is reasonable to look for a cubic $c$ in
$\rho^* Sym^3 T^*\M \subset Sym^3 T^* \widetilde{\M}$, satisfying
$\rfloor c = dp \circ i$.  But this is precisely where the {\it
Yukawa
cubic} of our CY family (also known as the Bryant-Griffiths
cubic, cf. [BG]) lives.  Abstractly, this is cup product,
followed by a ``rescaling'' by the gauge squared, and
integration, $$c: Sym^3H^1(T_X)
\stackrel{\mbox{cup}}{\rightarrow} H^3(\Lambda^3T_X) =
H^3(\omega^{-1}_X) \stackrel{s^2}{\rightarrow} H^3(\omega_X)
\stackrel{\int}{\rightarrow} \C.$$  Hodge theoretically, this is
the third iterate of the derivative of the period map (or of the
infinitesimal variation of Hodge structure): $$B_1 \circ B_2
\circ B_3: Sym^3H^1(T_X) \times H^{3,0} \to H^{0,3},$$ where the
$1$-dimensional spaces $H^{3,0}$ and $H^{0,3} \cong
(H^{3,0})^*$
are trivialized by means of $s$.  The equality $\rfloor c =
dp\circ i$
follows immediately from this interpretation, cf. [CGGH].

\vspace{0.2in}

     \underline{Step III}.  We still need to check that the
$2$-form we constructed is closed.  We could construct an affine
structure on $\widetilde{\M}$ and use Theorem 2.  Instead, we
verify the closedness directly, and obtain the affine structure
as corollary.  We do this by checking that the standard
symplectic form on $T^*\widetilde{\M}$ is invariant under
translation by the image of
$\lambda: {\cal H}_{3}(\tilde{{\cal X}}/\tilde{{\cal M}},\Z)
\hookrightarrow
(F^2{\cal H}^3)^*$, hence it descends to a symplectic form on
$\widetilde{\J}$.

     Since the question is local on $\widetilde{\M}$, we may
assume that the bundle of lattices $\HDI^3 \chi_*\Z$ on $\M$ is
trivial.  Let $X_0$ be the CY fiber over a base point $0 \in
\M$.  The choice of $\gamma_0 \in H_3 (X_0,\Z)$ determines a
family $\gamma = (\gamma_t)_{t \in \M}$, $\gamma_t \in
H_3(X_t, \Z)$, and these are taken by $\lambda$ to
$F^2H^3(X_t,\C)^*$.  By the identification $i$ of (5), we end
up
with a section $$\int_{\gamma} \in
\Gamma(\widetilde{\M},T^*\widetilde{\M}).$$ We need to show that
(the image of) this section is Lagrangian, with respect to the
standard symplectic form on $T^* \widetilde{\M}$.  Equivalently,
we
need to show that the $1$-form $\xi$ on $\widetilde{\M}$ which
corresponds to the section $\int_{\gamma}$ of
$T^*\widetilde{\M}$ is closed.

     Consider the function $g: \widetilde{\M} \to \C$ given by
$$g(X,s) := \int_{\gamma} s.$$ We claim $\xi = dg$.  Indeed, if
$a \in T_{(X,s)}\widetilde{\M}$ is the tangent vector to the
$1$-parameter family $(X_t,s_t)$ then we have, as in Step I:
$$s_t =
s_0 + ta \ \ \ \ (mod\ t^2).$$ and therefore
$$\begin{array}{lcl}<dg,a> & = & (\partial/\partial t)|_{t=0}
\ \ g(X_t,s_t)\\& = & (\partial/\partial t)|_{t=0}
\int_{\gamma}s_t\\& = & \int_{\gamma}a =
<\xi,a>.\end{array}$$\\
\begin{flushright}Q.E.D
\end{flushright}

\vspace{0.2in}

\noindent
\underline{Remark}.

     In order to obtain the Lagrangian structure, we did not need
to specify the ``level'' of the moduli space $\M$.  In order to
have an ACIHS, we need in addition to the Lagrangian structure
also a set of global Hamiltonians, i.e. a set of $1 + h^{2,1}$
independent functions on $\widetilde{\M}$.  For this, we choose
the level structure to include, at least, the choice of a
Lagrangian basis for the $H_3(X_t,\Z)$, say
$\gamma_0,\cdots,\gamma_{h^{2,1}}$.  As our Hamiltonians we then
take $$\tilde{t}_i := \int_{\gamma_i}s \ \ \ \ \ \ i =
0,\cdots,h^{2,1}.$$ We will see another interpretation for these
functions shortly.

\vspace{0.2in}

\proclaim{Lemma 2}.  The symplectic form $\sigma$ on
$\widetilde{\J}$ is exact. \par

\noindent
{\bf Proof}.

     For very general reasons, the symplectic form
$\tilde{\sigma}$ on $T^*\widetilde{\M}$ is exact: it equals
$d{\tilde{\tau}}$ , where $\tilde{\tau}$ is the action $1$-form
on $T^*\widetilde{\M}$.  (If $q_i$ are coordinates on
$\widetilde{\M}$ and $p_i$ the corresponding linear coordinates
on
the fibers, then $\tilde{\tau}$ is given by $\sum p_i dq_i$ and
$\tilde{\sigma}$ by $\sum dp_i \wedge dq_i$.)  To obtain
$\sigma$, we identified $T^*\widetilde{\M}$ with the vertical
bundle $\V$ of $\widetilde{\J}$ (whose fiber at $(X,s)$ is
$F^2H^3(X,\C))$, and observed that $\tilde{\sigma}$ was invariant
under translation by locally constant integral cycles $\gamma$,
hence it descended to $\sigma$ on $\widetilde{\J}$.

     A first guess for an antidifferential $\tau$ of $\sigma$
might be simply $\tilde{\tau}$ (pulled back by $i$, but we
suppress this), but $\tilde{\tau}$ is \underline{not} invariant
under translations: if the cycle $\gamma$ corresponds, as in
Step III of the proof of Theorem 3, to a $1$-form $\xi$ on
$\widetilde{\M}$, then translating by $\gamma$ changes
$\tilde{\tau}$ by addition of $\pi^*\xi$, where $\pi: \V \to
\widetilde{\M}$ is the projection.  To correct this discrepancy,
we consider the tautological function $f \in \Gamma(\O_\V)$ whose
value at a point $(X,s,v) \in \V$ (with $(X,s) \in
\widetilde{\M}$ and $v \in F^2H^3(X)^*$) is given by $$f(X,s,v)
:= v(s).\leqno(6)$$ Now $f$ is fiber-linear, so $df$ is constant
on fibers, and therefore translation by $\gamma$ changes $df$ by
$\pi^*$ of a $1$-form on $\widetilde{\M}$, and this $1$-form is
immediately seen to be $\xi$.  We conclude that $\tilde{\tau} -
df$ is a global $1$-form on $\V$ which is invariant under
translation by each $\gamma$, hence descends to a $1$-form
$\tau$ on $\widetilde{\J}$ which satisfies $d\tau = \sigma$, as
required.
\begin{flushright}Q.E.D
\end{flushright}

\vspace{0.2in}

\noindent
\underline{Remarks}.

     (1) Another way to see the exactness of $\sigma$ is to note
that it is the {\it symplectification} of a {\it contact
structure} $\kappa$ on $\J$.  In general, a contact manifold
$(\J,\kappa)$ determines an exact symplectic structure on the
natural $\C^*$-bundle $\widetilde{\J}$ over $\J$.  Conversely,
according to [AG], page 78, a symplectic structure $\sigma$ on a
manifold $\widetilde{\J}$ with a $\C^*$-action $\rho$ is the
symplectification of a contact structure on the quotient  $\J$
if and only if $\sigma$ is homogeneous of degree $1$ with respect
to $\rho$.  In our case, there are two independent $\C^*$-actions
on the total space of $\V \approx T^* \widetilde{\M}$: the
$\C^*$-action on $\widetilde{\M}$ lifts to an action $\rho'$ on
$T^*\widetilde{\M}$, and there is also the action $\rho''$ which
commutes with the projection to $\widetilde{\M}$ and is linear on
the fibers.  The symplectic form $\tilde{\sigma}$ is homogeneous
of weight $0$ with respect to $\rho'$ and of weight $1$ with
respect
to $\rho''$.  Neither of these actions descends to
$\widetilde{\J}$, but their product $\rho := \rho' \otimes
\rho''$ does.  The symplectic form $\sigma$ on $\widetilde{\J}$
is homogeneous of weight $1$ with respect to $\rho$, so it is the
symplectification of a contact structure on $\J$, and in
particular it is exact.

     In the notation of the proof, the actions are given, for $t
\in \C^*$, by :  $$\begin{array}{llcrrr}\rho' & : (X,s,v) &
\mapsto &
(X,& ts,& t^{-1}v)\\ \rho'' & : (X,s,v) & \mapsto & (X, & s, & t
v)\\ \rho & :
(X,s,v) & \mapsto & (X, & ts, & v).\end{array}$$ The
corresponding
vector fields on $\V$ are taken by $\tilde{\sigma}$ to the
$1$-forms $-df,\tilde{\tau}$, and $\tau = \tilde{\tau} - df$,
respectively.  Of these, only $\tau$ descends to a closed
$1$-form on $\widetilde{\J}$.  (Another description of this
contact
structure is in Remark (2) to Theorem 4.)

     (2) The tautological relative $3$-form $s$ on
$\widetilde{\X}$ over $\widetilde{\M}$ induces a relative
$1$-form on $\widetilde{\J}$ over $\widetilde{\M}$.
We note that our $1$-form $\tau$ on $\widetilde{\J}$ is a lift of
this tautological relative $1$-form.

\vspace{0.2in}

\proclaim{Lemma 3}.  The Hamiltonians $\tilde{t}_i =
\int_{\gamma_i}s$ (cf. remark after Theorem 3) are equal (up to
sign) to the action variables for $\widetilde{\J}$, given by (1).
\par

\vspace{0.2in}

\noindent
{\bf Proof}.

     The $3$-cycle $\gamma_i$ in $X$ determines a $1$-cycle, say
$\Gamma_i$, in $J(X)$.  The action variables (1) are given as
$\int_{\Gamma_i}\tau$.  By the proof of Lemma 2, $\tau$ pulls
back to $\tilde{\tau} - df$, where $f$ is given by (6),
and $\tilde{\tau}$ vanishes on fibers of $\pi$.  Therefore:
$$\int_{\Gamma_i} \tau = -\int_{\Gamma_i} df = -\int_{\gamma_i}
s.$$
\begin{flushright}Q.E.D
\end{flushright}

     Let $\A$ be the group of $1$-cycles in $X$ (i.e. the free
abelian group generated by all curves), and $\A_0$ the subgroup of
$1$-cycles homologous to $0$.  There is a well-defined {\it
Abel-Jacobi map}: $$\nu: \A_0 \to J(X).$$ If the cycle $C$
is the boundary of a $3$-chain $\Gamma$, then $\nu(C)$ is
the image in $J(X)$ of $\int_\Gamma$, considered as an
element of $H^3(X,\C)^*$.  The image is independent of the choice
of $\Gamma$.  (This is the analogue for $1$-cycles of the natural
map from divisors homologous to $0$
to $Pic^0(X)$ .)

     Given a family $\chi: \X \to \M$ of threefolds and a
relative
$1$-cycle ${\cal C} \to \M$ in $\A_0(\X/\M)$ (i.e. for each $t
\in \M$, $C_t$ is a $1$-cycle homologous to $0$ in $X_t$), there
is a corresponding section $$\nu({\cal C}) \in
H^0(\M,\J(\X/\M))$$ whose value at $t \in \M$ is $\nu(C_t)$
in $J(X_t)$.  This is called the {\it normal function} of
the relative cycle ${\cal C}$.

     Abstractly, a normal function is defined as a section
$\nu$ of $\J$ whose local liftings to sections
$\overline{\nu}$ of $\HDI^3\chi_*\C = {\cal H}^3$ satisfy the
differential equation $$\nabla \overline{\nu} \in F^1{\cal
H}^3,$$
or equivalently $$(\nabla \overline{\nu},s) = 0,\leqno(7)$$
where $\nabla$ is the Gauss-Manin connection on ${\cal H}^3$.
This
condition is independent of the choice of $\overline{\nu}$,
by Griffiths transversality.  It is automatically satisfied by
the Abel-Jacobi image of a family of relative $1$-cycles, cf.
[G2].

     More generally, we can consider multivalued $1$-cycles and
normal functions: a map $B \to \M$ gives a pullback family of
CYs over $B$, and we can consider a family of relative,
null-homologous $1$-cycles in this pullback; this determines a
normal
function over $B$, and we can consider its image (or that of an
abstract normal function over $B$) in $\J$.  We refer
to such subvarieties (no longer sections) of $\J$ as
multivalued normal functions.

\vspace{0.2in}

\proclaim{Theorem 4}.  For the symplectic structure on
$\widetilde{\J}$ constructed in Theorem 3, the image
$\nu(\A_0(\widetilde{\X}/\widetilde{\M}))$ is isotropic.  In
other words, the normal function of any (multivalued) family of
null-homologous $1$-cycles on Calabi-Yau threefolds gives an
isotropic subvariety of $\widetilde{\J}$.\par

\vspace{0.2in}

\noindent
{\bf Proof}.

     This is a variation on Step III in the proof of Theorem 3.
Say we are given a family ${\cal C} \to \B$ of null-homologous
$1$-cycles on gauged Calabi-Yaus; this means that there is a map
$t: \B \to \widetilde{\M}$ such that for $b \in \B$, the fiber
$C_b$ is a null-homologous $1$-cycle in $X_{t(b)}$.  The
Abel-Jacobi map $\nu: \B \to \widetilde{\J}$ is then a lift of
$t$, and we claim that its image is isotropic.  Locally in $\B$
we choose $3$-chains $\Gamma_b$ such that $\partial \Gamma_b =
C_b$; this determines local lifts $\overline{\nu}: \B \to
{\cal H}^3$ and $\tilde{\nu}: \B \to T^*\widetilde{\M}$ of
$\nu$.  Let $\tilde{\tau}$ be the action $1$-form on
$T^*\widetilde{\M}$, so $d\tilde{\tau}$ is the symplectic form.
We need to show that $\xi:= \tilde{\nu}^*\tilde{\tau}$ is a
closed $1$-form on $\B$.  At $b \in \B$, $\xi$ is given by
$\int_{\Gamma_b}$.  The new feature here is that instead of
$3$-cycles $\gamma_b \in H_3(X_b,\Z)$ we have the $3$-chains
$\Gamma_b$, with variable boundary $C_b$.

     As before, we think (locally in $\widetilde{\M}$) of $X$ as
being a fixed $C^\infty$ manifold, with variable complex
structure $\overline{\partial}_t$ and $3$-form $s_t$, $t \in
\widetilde{\M}$.  We consider the function defined locally on
$\B$: $$\begin{array}{lll}g & \colon & \B \longrightarrow
\C\\g(b) & := &
\int_{\Gamma_b} s_{t(b)};\end{array}$$ we claim that $\xi = dg$.
This time, both the integrand and the chain in $\int_{\Gamma_b}
s_{t(b)}$ depend on $b$.  Let $v$ denote an appropriate normal
vector to (the support of) $C_b$ along $\Gamma_b$.  (For details,
see [Gr2]).  We then have: $$\frac{\partial}{\partial b} g(b) =
\int_{\Gamma_b} \frac{\partial s_{t(b)}}{\partial b} + \int_{C_b}
v \rfloor s_{t(b)}.\leqno(8)$$ But $s_{t(b)}$ is of type $(3,0)$
with
respect to the complex structure $\overline{\partial}_{t(b)}$, so
the contraction $v \rfloor s_{t(b)}$ is of type $(2,0)$,
regardless of
the type of $v$.  On the other hand, the (Poincar\'e dual of) $C_b$
is of type $(2,2)$ with respect to this complex structure, so the
second term vanishes identically.  As in the proof of Theorem 3,
we then conclude that $dg = \xi$, as required.\\
\begin{flushright}Q.E.D
\end{flushright}

\vspace{0.2in}

\noindent
\underline{Remarks}.  The proof shows that Theorem 4 can be
strengthened in two directions:

     (1) Let $\tau = \tilde{\tau} - df$ be the $1$-form on
$\tilde{\J}$ constructed in Lemma 2, satisfying $d\tau = \sigma$.
(The function $f$ is given by formula (6).)  Then each normal
function $\nu$ gives an integral manifold of $\tau$, i.e.
$\nu^*\tau = 0$.  (This of course implies $\nu^*\sigma
= 0$.)  Indeed, in the notation of the proof we have $\nu^*
\tilde{\tau} = \xi$ and $\nu^*f = g$, so what we show is
exactly the vanishing of $\nu^*\tau = \xi - dg$.  If we
consider instead normal functions for $\J$ over $\M$ (suppressing
the gauge), we see that they are integral manifolds for the
contact structure on $\J$ given in Remark 1 after Lemma 2.

     (2) The result of Theorem 4 holds for any normal function,
not only for those coming from cycles.  In fact, the condition
$\nu^*\tau = 0$ is {\it equivalent} to Griffiths'
differential equation for a normal function.  In this general
setting we define the function $g$ on $\B$ as the inner product
$$g := (\overline{\nu},s)$$ of the sections
$\overline{\nu}$ (= the lift of $\nu$) and $s$ (= the
$(3,0)$-form).  Its value depends only on the partial lift
$\tilde{\nu}: \B \to T^*\widetilde{\M}$; we have $g =
\tilde{\nu}^*f$.  The crucial formula (8) becomes: $$dg =
(\overline{\nu},\nabla s) + (\nabla
\overline{\nu},s).$$ So we have $\nu^*\tau =
\tilde{\nu}^* \tilde{\tau} - dg = -(\nabla
\overline{\nu},s)$, and the vanishing of $\nu^*\tau$ is
indeed equivalent to the differential equation (7) of a normal
function.

     We can define the contact structure on $\J$, and hence the
symplectic structure on $\widetilde{\J}$, in terms of normal
functions: The contact distribution is given at a point
$(X_0,\nu_0) \in \J$ as the subspace of
$T_{(X_0,\nu_0)}\J$ spanned by differentials of all normal
functions $\nu$ defined near $X_0$.  If we choose one such
$\nu$, we obtain a splitting $$T_{(X_0,\nu_0)} \J
\approx F^2H^3(X_0)^* \oplus T_{X_0}\M,$$ in terms of which the
contact distribution is $$H^{2,1}(X_0)^* \oplus T_{X_0}\M,$$ and
this subspace is independent of the choice of $\nu$, by (7).

     (3) An example of a ``multi-valued'' normal function is
given by the inclusion $\nu: \B \hookrightarrow
\widetilde{\J}$ of the total space $\B$ of a family of abelian
subvarieties, $\B \to S$ (where $S$ comes with a map $i: S \to
\widetilde{\M}$, and the fiber $B_s$ is an abelian subvariety of
the torus $J(X_{i(s)})$).  We conclude:

\vspace{0.2in}

\proclaim{Corollary}.  The total space of any family of abelian
subvarieties of $\widetilde{\J}$ is isotropic.\par

     The generalized Hodge conjecture says that $\nu(\B)$ is
actually contained in
$\nu(\A_0(\widetilde{\X}/\widetilde{\M})$), but this is not
known.  At first sight it also seems plausible that a maximal
$\nu(\B)$, as well as each component of
$\nu(\A_0(\widetilde{\X}/\widetilde{\M}))$, should actually
be {\it Lagrangian}; there are, however, counterexamples to this.

     The intermediate Jacobian is the connected component of the
origin in a larger group $D(X)$, the Deligne cohomology group of
$X$, cf. [EZ].  This fits in the exact sequence $$0 \to J(X) \to
D(X) \stackrel{\varphi}{\rightarrow} H^{2,2} (X,\Z) \to 0,$$
where $H^{2,2}(X,\Z) := H^{2,2}(X,\C) \cap H^4(X,\Z)$.  For a
Calabi-Yau, this is the same as $H^4(X,\Z)$.  (The relation of
$D$ to $J$ is analogous to that of $Pic$ to $Pic^0$.)  The
Abel-Jacobi map $\A_0(X) \to J(X)$ extends to an Abel-Jacobi
(aka cycle-class) map, $$\nu: \A(X) \to D(X),$$ such that
$\varphi \circ \nu: \A(X) \to H^{2,2}(X,\Z)$ sends a cycle
to its fundamental class.  Everything works well for families: a
family $\X \to \M$ determines a relative Deligne group $\D(\X/M)$
whose fiber over $t \in \M$ is $D(X_t)$, and there is an
Abel-Jacobi map $$\nu: \A(\X/\M) \to \D(\X/\M)$$ which
fiber-by-fiber agrees with the previous $\nu$.

\vspace{0.2in}

\proclaim{Theorem 5}.  Let $\D = \D(\X/\M)$ be the relative
Deligne group of a complete family of CY threefolds, let
$\widetilde{\D} = \D(\widetilde{\X}/\widetilde{\M})$ be the
relative Deligne group of the gauged family, and let $\J,
\widetilde{\J}$ be the corresponding relative Jacobians.  Then
there is a natural contact structure $\kappa$ on $\D$ with
symplectification $\sigma = d\tau$ on $\widetilde{\D}$ with the
following properties:
\begin{list}{{\rm(\alph{bean})}}{\usecounter{bean}}
\item $\sigma, \tau$ and $\kappa$ restrict to the previously
constructed structures on $\widetilde{\J}$ and $\J$.
\item The fibration $\widetilde{\D} \to \widetilde{\M}$ is
Lagrangian.
\item All normal functions
$\nu(\A(\widetilde{\X}/\widetilde{\M}))$ are integral
manifolds for $\tau$ (hence isotropic for $\sigma$), and normal
functions $\nu(\A(\X/\M))$ are integral manifolds for
$\kappa$.
\end{list}

\vspace{0.2in}

\noindent
{\bf Proof}.

     Each component of $\D$ (respectively $\widetilde{\D}$) is
locally isomorphic to $\J$ (respectively $\widetilde{\J}$.  In
order to extend $\sigma,\tau,\kappa$ we need a collection of
local sections of $\D$ (and $\widetilde{\D}$) with the
properties:
\begin{list}{{\rm(\arabic{bean})}}{\usecounter{bean}}
\item For each $t \in \M$ and each component of $D(X_t)$, there
is a local section defined near $t$ and passing through this
component.
\item The difference of any two local sections in the same
component is an integral manifold for $\sigma,\tau, \kappa$.
\end{list}

     Indeed, such a collection of sections determines a
collection of local isomorphisms of components of $\D$
(respectively $\widetilde{\D}$) with $\J$ (respectively
$\widetilde{\J}$), so we define $\sigma,\tau,\kappa$ as pullbacks
of the corresponding structures on $\widetilde{\J}$ and $\J$.
This is well-defined by (2), and covers $\widetilde{\D}$ and
$\D$, by (1).  Now such a collection is given, in our case, by
the normal functions of all $1$-cycles, $\nu(\A(\X/\M))$:
property (1) follows from the Lefschetz theorem on
$(1,1)$-classes, and property (2) amounts to our Theorem 4 and
the
remarks following it.\\
\begin{flushright}Q.E.D
\end{flushright}

     As an immediate corollary, we get the following result of
Griffiths [G1].  (There it is proved only for quintic threefolds,
but the argument works in general, once we know the
unobstructedness of moduli [Bo,Ti,To].)

\vspace{0.2in}

\proclaim{Corollary}.  On a generic CY, algebraic equivalence
of $1$-cycles implies $AJ$-equivalence; in other words, the
Abel-Jacobi map is constant on any connected family of
$1$-cycles.\par

\vspace{0.2in}

\noindent
{\bf Proof}.

     The $AJ$ image of the family is isotropic, by Theorem 5, and
dominates $\widetilde{\M}$, hence it is Lagrangian and meets the
generic fiber discretely, so the $AJ$ map is locally constant on
the connected family, hence constant.\\
\begin{flushright}Q.E.D
\end{flushright}

\vspace{0.2in}

\noindent
{\bf Note}.  These results have analogues for Calabi-Yaus of odd
dimensions $n = 2k+1\geq 3$:\ The universal $k^{th}$ intermediate
(or middle-dimensional) Jacobian
$\widetilde{\J}^k(\widetilde{\X}/\widetilde{\M})$, as well as the
corresponding relative Deligne cohomology, have natural
quasi-symplectic structures, i.e. closed $2$-forms which in
general
will be degenerate.  The fibers of the projections of
$\widetilde{\J}$ or $\widetilde{\D}$ to $\widetilde{\M}$ are
Lagrangian, i.e. maximal isotropic with respect to the $2$-form
(but of dimension greater than $\dim \widetilde{\M}$, because of
the degeneracy).  All normal functions, or Abel-Jacobi images of
families of $k$-cycles, are isotropic.  For the details, see
[DM].

\vspace{0.3in}

\S 3. {\bf Mirrors}.

     To a given family $\X \to \M$ of Calabi-Yaus, there are two
ways of associating an ($N=2$, $c=9$) (super) {\it conformal field
theory} (CFT).  Physicists refer to these as the $A$-model and
$B$-model theories.  A key quantity in the CFT is its {\it
partition function}, which is a cubic form on a vector space $V$,
depending on some auxiliary parameters $q$.

     In the $A$-model, one takes $V := H^{1,1}(X) = H^2(X,\C)$,
where $X$ is a general CY in the given family.  On it there is
the cubic form Topo, giving the topological cup product.  For
each $k \in H_2(X,\Z) \approx H^2(X,\Z)^*$ (we replace
$\Z$-homology or cohomology by its image in the
${\bf Q}$-(co)homology, i.e. we ignore torsion classes), we think
of $k^3$
as a cubic form on $V$.  The partition function is then: $$Topo +
\sum_{k,d} n_k k^3 q^{dk}.\leqno(9)$$ In the summation, $k$ runs
over $H^2(X,\Z)^*$ and $d$ over positive integers.  The
coefficient $n_k$ is the ``number of rational curves'' of
fundamental class $k$ in $X$.  This number is conjectured, but
not known, to be finite.  The multi-parameter $q$ is defined by
$q := e^{2\pi it}$, where $t$ is a set of linear coordinates on
$H^2(X,\C)$.  Intrinsically, we define $q^{dk}$ as the function
given at $t \in H^2(X,\C)$ by: $$q^{dk} := e^{2\pi
id<t,k>}.\leqno(10)$$ The sum is expected to converge for $t$ in
$K$, the complexified K\"{a}hler cone of $X$, i.e. when $Im(t)
\in H^2(X,\R)$ is in the K\"{a}hler cone.

     In the $B$-model, $V$ is $H^{2,1}(X)$, and the partition
function is a normalized version of the Yukawa cubic, so its
coefficients are functions on $\M$.  The Yukawa cubic is defined
only up to a scalar factor which is determined by a choice of
holomorphic gauge $s \in H^{3,0}(X)$.  In order to get the
coefficients in the partition function to be functions on $\M$
(rather than $\widetilde{\M}$), we need a section of
$\widetilde{\M} \to \M$.  According to [Mo], one chooses a {\it
maximally unipotent degeneration} of $X$; this determines a
vanishing cycle $$\gamma_0 \in H_3(X,\Z),$$ and $s$ is
normalized by the condition $$\int_{\gamma_0} s_0 = 1.$$

     Physicists believe [Gep] that either the $A$-model or
$B$-model constructions set up bijections between CFTs and CY
families.  This, together with the existence of an internal
symmetry of CFT which interchanges $A$ and $B$ variables, led
[Di,LVW] to the {\it mirror conjecture}, which says that each CY
family $\X \to \M$, perhaps with some extra data, determines
another family $\X' \to \M'$ such that the $A$-model CFT of the
former is isomorphic to the $B$-model CFT of the latter, and vice
versa.

     It is not clear (to us) what the precise mathematical
formulation of this conjecture should be.  What is clear, though,
is that the conjecture involves the existence of an isomorphism,
called the {\it mirror map}. $$t: \M \to K'\leqno(11)$$ where
$\M$ is (an open subset of ?) the moduli space of Calabi-Yaus
(with extra structure ?) in the first family, and $K'$ is the
complexified K\"{a}hler cone of an $X'$ in the second family.
This should satisfy:
\begin{list}{{\rm(\alph{bean})}}{\usecounter{bean}}
\item $dt: H^1(T_X) \stackrel{\sim}{\rightarrow} H^{1,1}(X')$ is
an isomorphism.  (In particular, the dimensions agree, so that
the Hodge diamonds of $X,X'$ are obtained from each other by a
$90^{\circ}$ turn).
\item The partition functions correspond through $t$, i.e.:
$${\mbox Yukawa}_{(X,s_0)} = (Topo_{X'} +
\sum_{\stackrel{k \in H_2(X',\Z)}{d \in \Z_+}} n_k(X')k^3q^{dk})
\circ
dt,$$ where both sides are interpreted as sections of
$Sym^3T^*_{\M}$.
\end{list}

\vspace{0.2in}

     In [CdOGP], this conjectural equality of partition functions
was used to obtain some spectacular predictions for numbers of
rational curves of arbitrary degrees in quintic hypersurfaces in
${\bf P}^4$.  These predictions have since been extended to many
other CY families.  We refer to [Mo] for more details of the
conjecture and the evidence for it.

     In the remainder of this section we wish to speculate on the
role which the ACIHS on relative Deligne cohomology
$\widetilde{\D}$ might play in the mirror story.  A couple of
relationships are evident:
\begin{list}{{\rm(\roman{bean})}}{\usecounter{bean}}
\item Both of the relevant dimensions, $h^{1,1}$ and $h^{2,1}$,
occur in $\widetilde{\D}$, as the rank of the discrete part and
(one less than) the dimension of the continuous part (of either
base or fibers), respectively.  This suggests that the relation
between the two ACIHS attached to mirror-symmetric CY families
might be expressible directly as some sort of Fourier transform.
\item The mirror map, $t$ of (11), as constructed in [CdOGP],
[Mo], etc., is given by $$t_i = \tilde{t}_i/\tilde{t}_0, \ \ \ \
i = 1,\cdots,h^{2,1},$$ where $\tilde{t}_i$ are the action
coordinates for the ACIHS, for an appropriate choice of
Lagrangian basis $\gamma_0,\cdots,\gamma_{h^{2,1}}$ in terms of
a maximally unipotent degeneration.  (This follows immediately
from the construction in [Mo] and our Lemma 3.)
\item We also recall from [Mo] that mirror symmetry is expected
to hold for ``most'' Calabi-Yaus, but probably not all; a
notorious collection of probable exceptions is provided by
Calabi-Yaus which are rigid, i.e. $h^{2,1} = 0$.  Now a somewhat
analogous situation occurs for the Torelli and Schottky problems,
of recovering a variety or family of varieties from its Hodge
structure: this is known to be true, at least generically, for
many classes of varieties, but various counterexamples exist.  It
is conceivable that all failures of mirror symmetry could be
explained by failures of Schottky or Torelli results, i.e. that
there is still some mirror data, but it does not come from a
geometric object, namely a CY family.
\end{list}

     The above observations, and some others which will be
discussed below, suggest that we might try to split the mirror
conjecture roughly into three parts:
\begin{list}{{\rm(\arabic{bean})}}{\usecounter{bean}}
\item Axiomatize the integrable system associated to a CY family,
by listing a collection of cohomological and symplectic data
which can be attached to each family, and which satisfies certain
axioms.
\item Describe a formal mirror transform on abstract systems,
which interchanges the continuous and discrete parts.
\item Study the validity of Torelli- and Schottky-type results,
which say that a CY-family can be uniquely recovered from its
system, and describe which abstract systems actually arise from
geometry.
\end{list}

     A weaker version of this last step, likely to hold in
greater generality, would be:

(3') Extract from an abstract system enough data to reconstruct
either the $A$- or $B$-model partition functions.  (Ideally, we
should be able to extract both partition functions from the
system, and the formal transform of step (2) should interchange
them.)

     An abstract system should consist at least of the following
data:
\begin{list}{{\rm(\roman{bean})}}{\usecounter{bean}}
\item An analytically completely integrable Hamiltonian system
$\widetilde{\D} \to \widetilde{\M}$, which is an extension of
another ACIHS $\widetilde{\J} \to \widetilde{\M}$ by a lattice
$H^{1,1}(\Z)$.
\item An indefinite (``Lorentzian'') polarization on the torus
fibers of $\widetilde{\J} \to \widetilde{\M}$ and a (symmetric,
cubic, $\Z$-valued) intersection product on $H^{1,1}(\Z)$.
\item A $\C^*$-action on $\widetilde{\D}$ lifting an action of
$\C^*$ on $\widetilde{\M}$.
\item A $1$-form $\tau$ on $\widetilde{\D}$, homogeneous of
weight $1$ with respect to the $\C^*$-action, satisfying $d\tau =
\sigma$ (= the symplectic form on $\widetilde{\D}$).
\end{list}

     Some further pieces of data may also need to be specified,
such as the choice of a maximally unipotent degeneration.  The
above data lets us recover the variation of weight-$3$ Hodge
structures over $\widetilde{\M}$, and the family of all normal
functions: these are the integral manifolds for $\tau$.  A key
question which we cannot answer is whether the collection of
normal functions coming from $1$-cycles is determined by the
above data, or needs to be specified explicitly.

     These data certainly seem rich enough that one expects to
recover from them much of the underlying CY geometry.  The hope,
which of course we are nowhere near realizing, is that when
these data satisfy some appropriate axioms they can be used to
generate a mirror system of the same type.  If this is so, one
must be able to extract both the $A$- and $B$-model data (i.e.
the two partition functions) from a given system.  Let us discuss
the various ingredients.

     For the $B$-model, we need two essential pieces: the Yukawa
cubic $c$ and the affine structure on $\widetilde{\M}$.  Now the
Yukawa cubic is determined by the symplectic form on
$\widetilde{\J}$ and the polarization on the fibers, by Theorem
1.
The affine structure used to convert $c$ to a cubic on the fixed
vector space $V$ is determined by the coordinates $\tilde{t}_i$,
which by Lemma 3 are the action variables, so they too are
determined by $\widetilde{\J}$.  This of course depends on the
discrete choice of a Lagrangian basis for $H_3(X,\Z)$; we do not
see how to avoid this choice, so it must either be built into the
axioms, or be deducible somehow from knowing the collection of
all normal functions (and in particular, their monodromy).

     For the $A$-model it seems plausible that the relevant data
can be obtained, but this leads us to Hodge-theoretic questions
which we cannot (yet?) resolve.  We need the topological
intersection form Topo, which is built into the axioms.  The
essential problem is to recover the coefficients $n_k$.  Each
rational curve on the generic CY in the family gives a (finitely
multivalued) normal function in $\widetilde{\D}$, and each
homologous pair gives one in $\widetilde{\J}$.  The problem is
that not all normal functions arise in this way.  We are led to a
series of problems.
\begin{list}{{\rm(\roman{bean})}}{\usecounter{bean}}
\item Which normal functions in $\widetilde{\D}$ arise from
effective curves?
\item Which of these arise from irreducible curves?
\item Which of these arise from rational curves?
\end{list}

     We note that a given normal function $\nu$ might arise
from many different curves.  The simplest example is that of
complete intersection curves: a given family of such curves is
parametrized by a (usually positive dimensional) projective space
on each CY in the family, and all such curves have the same
Abel-Jacobi image since they are rationally equivalent.  In fact,
the
Corollary to Theorem 5 says that all curves in any connected
family on the generic CY give rise to the same $\nu$.

     The above remark suggests that we need not only decide which
normal functions $\nu$ come from rational curves, but to
determine a multiplicity or weight of $\nu$ accounting for
all its realizations in terms of rational curves.  Now this seems
rather special; one can just as well hope for a weight function
$w_{\nu}$ counting all curves, of all genera, corresponding
to $\nu$, through an expansion such as $w_\nu = \sum
N_{\nu,g} \lambda^{g-1}$, where $\lambda$ is a dummy
variable and the coefficients $N_{\nu,g}$ account for the
curves.

     A master partition function $\F$, accounting for curves of
all genera on a CY, has recently been proposed in [BCOV].  It can
be written:
$$\begin{array}{ccl}\F & = & \sum F_g \lambda^{g-1}\\& = &
\sum N_{k,g}q^k \lambda^{g-1}\\& = & \sum w_k q^k,\\
\end{array}$$ where the $q^k$ are as in (10), the $F_g$ are the
genus-$g$ partition functions (they depend on the $q^k$), and the
coefficients $N_{k,g}$ are certain combinations of the
``numbers'' $n_{k,g}$ of curves of genus $g$ and fundamental
class $k$.  [A rigid non-singular curve of class $k$ and genus
$g$ contributes a term of $1 \cdot \lambda^{g-1}$, but also makes
contributions to higher-genus terms, because of multiple
(branched) coverings.  The situation is more complicated for
continuous families of curves, where the basic contribution
involves the Euler characteristic of the base, but there are
other terms corresponding, e.g., to singular curves in the
family.] This $\F$ is a function of the variables $q$ and
$\lambda$.  Via the mirror correspondence (11), it is carried to
a function on $\widetilde{\M}$, whose leading term is the Yukawa
coupling.  The genus-$0$ partition function $F_0$ is recovered as
a residue of $\F$ in the $\lambda$ direction.

     In order to recover the partition function $F_0$, or $\F$
for that matter, we need a refinement of the BCOV coefficients:
to each Lagrangian normal function $\nu$ we want to
associate a weight function $$w_\nu = \sum N_{\nu,g}
\lambda^{g-1}$$ such that for each nef class $k$ we have $$w_k =
\sum_{\{\nu|\varphi(\nu) = k\}}w_\nu,$$ or $$N_{k,g} =
\sum_{\{\nu|\varphi(\nu)=k\}} N_{\nu,g},$$ or equivalently:
$$\F = \sum N_{\nu,g} q^{\varphi(\nu)} \lambda^{g-
1}.$$ The point is that we want a formula for the $w_\nu$ in
terms of $\nu$ alone, without knowledge of the actual curves
on $X$ whose Abel-Jacobi image is $\nu$.

     We do not know whether such a formula exists.  But if it
does, it is likely to involve only the {\it infinitesimal
invariant} $\delta \nu$, associated to a normal function by
Griffiths [G2] and improved by Green [Gr1].  This invariant lives
in a certain Koszul cohomology group.  It vanishes if and only if
$\nu$ is locally constant, i.e. can be lifted locally to a
constant section $\overline{\nu}$ of the Gauss-Manin local
system $\HDI^3 \widetilde{\chi}_* \C$ of the family
$\widetilde{\chi}:
\widetilde{\X} \to \widetilde{\M}$.

     The reason is that the sum, over all $\nu$ in
$\varphi^{-1}(k)$, is too big.  For example, when $k=0$, any
torsion element in $J(X)$ (or in $H^3(X,{\bf Q}/\Z))$ extends to
a finitely valued Lagrangian multisection of $\widetilde{\J}$,
but there is no reason to expect any contribution of these
torsion sections to the curve count.  By continuity, one expects
no contribution from any locally constant normal function, so
$w_\nu$ should depend only on $\delta\nu$.

     In our case, of normal functions on CY threefolds, we can describe the
Koszul cohomology group in which $\delta\nu$ lives explicitly in terms
of the Yukawa cubic.  In general, $\delta\nu$ is obtained by
choosing a lift $\overline{\nu}$, differentiating it with
respect to the Gauss-Manin connection $\nabla$, and taking the
part independent of choices.  In our case, $\nabla \overline{\nu}$ lives in
$$Hom(T_{X,s} \widetilde{\M}, H^3/F^2H^3)$$ which is canonically
identified with $\otimes^2(F^2)^*$.  The differential equation
(7) satisfied by any normal function asserts that
$\overline{\nu}$ lands in
$$Hom(T_{X,s}\widetilde{\M},F^1H^3/F^2H^3)$$ which is $(F^2)^*
\otimes (H^{2,1})^*$.
By Remark (2) to Theorem 4, any normal function is Lagrangian.
This means that $\nabla \overline{\nu}$ is symmetric,
so we are down to $$Sym^2(H^{2,1})^*,$$ or the space of quadrics
on $V = H^{2,1}$.  The ambiguity arising from a change of the lift
$\overline{\nu}$ amounts to those quadrics which are partial
derivatives, in some direction, of the Yukawa cubic $c \in
Sym^3(H^{2,1})^*$.
Summarizing, the Koszul complex in our case can be identified
with the complex
\[ (H^{2,1}) \stackrel{Yukawa}{\rightarrow} \otimes^{2} (H^{2,1})^*
\rightarrow \stackrel{2}{\wedge}(H^{2,1})^*.
\]
We conclude that $\delta\nu$ is in $R^2$,
where $R^d$ is the $d^{th}$ graded piece of the Jacobian ring of
$c: R^{\bullet} := S^{\bullet}/\J^{\bullet}$, where $S^{\bullet}
= Sym^{\bullet}(V^*)$
and $\J^{\bullet} \subset S^{\bullet}$ is the Jacobian ideal,
generated by the
partials of $c$.  The problem of reconstructing the curves on $X$
with given normal function $\nu$ from the invariant
$\delta\nu$ may therefore have some resemblance to the
variational Torelli problem, where certain parts of $R$ are
given and one needs to recover the underlying variety, as in [D].
It is known how to extract quite a bit of information from the
invariants $\delta\nu$ in some situations, e.g. [G2], [Gr1],
[CP].  The question is whether, on a generic Calabi-Yau threefold
$X$, one can recover the collection of all curves from their
$\delta\nu$ invariants.

\Bigskip

\centerline{\bf REFERENCES}

\medskip

\def\bib#1{\noindent\hbox to50pt{[#1]\hfil}\hang}
\def\bibline#1{\bib{#1}\vrule height.1pt width0.75in depth.1pt
\/,}

\vskip 25pt
\parindent=50pt
\frenchspacing

\bib{AG}V.I. Arnol'd and A.B. Givental': {\it Symplectic
Geometry}. In: Arnol'd, V.I., Novikov, S.P. (eds.) Dynamical
systems IV (EMS, vol. 4), Springer (1988), 1-136.

\bib{AHH}M.R. Adams, J. Harnad, J. Hurtubise: {\it Isospectral
Hamiltonian flows in finite and infinite dimensions, II.
Integration of flows}. Commun. Math. Phys. {\bf 134} (1990),
555-585.

\bib{B}A. Beauville: {\it Jacobiennes des courbes spectrales et
syst\`{e}mes hamiltoniens compl\`{e}tement int\'{e}grables}. Acta
Math. {\bf 164} (1990), 211-235.

\bib{BCOV}M. Bershadsky, S. Cecotti, H. Ooguri and C. Vafa:
{\it Holomorphic Anomalies in Topological Field Theories}, with
an
appendix by S. Katz, Harvard University preprint HUTP-93/A008,
RIMS-915.

\bib{BG}R. Bryant and P. Griffiths: {\it Some observations on the
infinitesimal period relations for regular threefolds with
trivial canonical bundle}, in Arithmetic and Geometry, Papers
dedicated to I.R. Shafarevich, vol. II, Birkh\"auser, Boston
(1983), 77-102.

\bib{Bn}F. Bottacin: Thesis, Orsay, 1992.

\bib{Bo}F. Bogomolov: {\it Hamiltonian K\"{a}hler manifolds},
Soviet Math Dokl. {\bf 19} (1978), 1462-1465.

\bib{CdOGP}P. Candelas, X.C. de la Ossa, P.S. Green, and L.
Parkes: {\it A pair of Calabi-Yau manifolds as an exactly soluble
superconformal theory}, Phys. Lett. {\bf B 258} (1991), 118-126;
Nuclear Phys. {\bf B 359} (1991), 21-74.

\bib{CGGH}J. Carlson, M. Green, P. Griffiths, and J. Harris: {\it
Infinitesimal variations of Hodge structure, I}, Compositio Math.
{\bf 50} (1983), 109-205.

\bib{CP}A. Collino, G. P. Pirola: {\it
Griffiths' infinitesimal invariant for a curve in its Jacobian},
Preprint.

\bib{D}R. Donagi: {\it Generic Torelli for projective
hypersurfaces}, Compositio Math. {\bf 50} (1983), 325-353.

\bib{DM}R. Donagi and E. Markman: {\it Spectral curves,
algebraically completely integrable Hamiltonian systems, and
moduli of bundles}, CIME lecture notes, 1993. To appear in LNM.

\bib{Di}L.J. Dixon: {\it Some world-sheet properties of
superstring compactifications, on orbifolds and otherwise},
Superstrings, Unified Theories, and Cosmology (1987) (G. Furlan
et al., eds.), World Scientific, Singapore, New Jersey, Hong Kong
(1988),  67-126.

\bib{DKN}B.A. Dubrovin, I.M. Krichever and S.P. Novikov: {\it
Integrable systems I}, In: Arnol'd, V.I., Novikov, S.P. (eds.)
Dynamical
systems IV (EMS, {\bf vol. 4}), Springer (1988), 173-280.

\bib{EZ}F. El Zein and S. Zucker: Extendability of normal
functions
associated to algebraic cycles. In {\it Topics in transcendental
Algebraic Geometry}, ed. by P. Griffiths, Annals of Mathematics
Studies {\bf 106} Princeton (1984), 260-288.

\bib{G1}P.A. Griffiths: {\it On the periods of certain rational
integrals I,II}; Ann. of Math. {\bf 90} (1969), 460-541.

\bib{G2}P.A. Griffiths:  {\it Infinitesimal variations of Hodge
structure III: determinantal varieties and the infinitesimal
invariant of normal functions}; Compositio Math. {\bf 50} (1983),

267-324.

\bib{Gep}D. Gepner: {\it Exactly solvable string compactification
on manifolds of $SU(N)$ holonomy}, Phys. Lett. {\bf B 199}
(1987), 380-388.

\bib{Gr1}M. Green: {\it Griffiths infinitesimal invariant and the
Abel-Jacobi map}, J. Diff. Geometry {\bf 29} (1989), 545-555.

\bib{Gr2}M. Green: {\it Infinitesimal methods in Hodge theory},
CIME notes, June 93, to appear in LNM.

\bib{H}N.J. Hitchin: {\it Stable bundles and integrable systems},
Duke Math. J. {\bf 54} (1987),  91-114.

\bib{LVW}W. Lerche, C. Vafa and N.P. Warner: {\it Chiral rings
in $N = 2$ superconformal theories}, Nuclear Phys. {\bf B 324}
(1989), 427-474.

\bib{M1}E. Markman: {\it Spectral curves and integrable systems},
to appear in Comp. Math.

\bib{M2}E. Markman: In preparation.

\bib{Mo}D. Morrison: {\it Mirror Symmetry and Rational Curves on
Quintic Threefolds: A Guide for Mathematicians}, J. Am. Math.
Soc. {\bf 6} (1993), 223.

\bib{Ti}G. Tian: {\it Smoothness of the universal deformation
space of compact Calabi-Yau manifolds and its Petersson-Weil
metric}, in Mathematical Aspects of String Theory (S.T. Yau,
ed.), World Scientific, Singapore (1987), 629-646.

\bib{To}A. Todorov: {\it The Weil-Petersson geometry of the moduli
space of $SU(n \geq 3)$ (Calabi-Yau) manifolds, I}, Comm. Math.
Phys. {\bf 126} (1989), 325-346.

\bigskip
\noindent
e-mail addresses: donagi@math.upenn.edu, markman@math.lsa.umich.edu
\end{document}